\documentclass[a4paper,10pt]{amsart}

\usepackage{verbatim}
\usepackage{enumerate}
\usepackage{amsmath, amsthm, amssymb}
\usepackage{amsfonts}
\usepackage{mathrsfs}

\usepackage{hyperref}
\usepackage[active]{srcltx}

\newtheorem{thm}{Theorem}[section]
\newtheorem{lem}[thm]{Lemma}
\newtheorem{defn}[thm]{Definition}
\newtheorem{prop}[thm]{Proposition}
\newtheorem{cor}[thm]{Corollary}

\theoremstyle{remark}
\newtheorem{ex}[thm]{Example}
\newtheorem{rem}[thm]{Remark}

\newcommand{\ben} {\begin{enumerate}}
\newcommand{\beni} {\begin{enumerate}[\rm(i)]}
\newcommand{\een} {\end{enumerate}}
\newcommand{\beq} {\begin{equation}}
\newcommand{\eeq} {\end{equation}}
\newcommand{\beqw} {\begin{equation*}}
\newcommand{\eeqw} {\end{equation*}}

\renewcommand{\d}{\delta}
\newcommand{\g}{\gamma}

\renewcommand{\l}{\lambda}
\newcommand{\om}{\omega}
\newcommand{\Om}{\Omega}

\newcommand{\s}{^*}

\newcommand{\e}{\varepsilon}

\newcommand{\E}{\mathbb{E}}
\newcommand{\R}{\mathbb{R}}

\newcommand{\F}{\mathscr{F}}
\renewcommand{\P}{\mathbb{P}}
\renewcommand{\L}{\mathscr{L}}

\newcommand{\dD}{\mathbb{D}}

\newcommand{\D}{\rm{dom}}
\newcommand{\one}{\mathbf{1}}

\newcommand{\n}{\Vert}
\newcommand{\lb}{\langle}
\newcommand{\rb}{\rangle}

\newcommand{\bpr}{\begin{proof}}
\newcommand{\epr}{\end{proof}}
\newcommand{\blem}{\begin{lem}}
\newcommand{\elem}{\end{lem}}
\newcommand{\bdefn}{\begin{defn}}
\newcommand{\edefn}{\end{defn}}
\newcommand{\bthm}{\begin{thm}}
\newcommand{\ethm}{\end{thm}}
\newcommand{\bprop}{\begin{prop}}
\newcommand{\eprop}{\end{prop}}
\newcommand{\bcor}{\begin{cor}}
\newcommand{\ecor}{\end{cor}}
\newcommand{\bex}{\begin{ex}}
\newcommand{\eex}{\end{ex}}
\newcommand{\brem}{\begin{rem}}
\newcommand{\erem}{\end{rem}}

\renewcommand{\H}{\mathscr{H}}
\newcommand{\calL}{\mathscr{L}}
\newcommand{\calS}{\mathscr{S}(\Om)}
\newcommand{\calSE}{\mathscr{S}(\Om)\ot E}
\newcommand{\calSc}{\mathscr{S}_{\rm c}(\Om)}
\newcommand{\calSEc}{\mathscr{S}_{\rm c}(\Om)\ot E}
\newcommand{\calSEs}{\mathscr{S}(\Om)\ot E^*}

\newcommand{\ip}[1]{\langle {#1}\rangle}
\newcommand{\biggip}[1]{\bigg\langle {#1}\bigg\rangle}

\newcommand{\ot}{\otimes}

\begin{document}

\title{A Clark-Ocone formula in UMD Banach spaces}

\author{Jan Maas}
\address{Delft Institute of Applied Mathematics\\
Delft University of Technology\\ P.O. Box 5031\\ 2600 GA Delft\\The
Netherlands}
\email{J.Maas@tudelft.nl}

\author{Jan van Neerven}
\address{Delft Institute of Applied Mathematics\\
Delft University of Technology\\ P.O. Box 5031\\ 2600 GA Delft\\The
Netherlands}
\email{J.M.A.M.vanNeerven@tudelft.nl}

\date\today

\thanks{The authors are supported by VIDI subsidy 639.032.201
(JM, JvN) and VICI subsidy 639.033.604 (JvN)
of the Netherlands Organization for Scientific Research (NWO).
The first named author
acknowledges partial support by the ARC Discovery Grant DP0558539}

\begin{abstract}
Let $H$ be a separable real Hilbert space and let
$\mathbb{F}=(\F_t)_{t\in [0,T]}$ be the augmented filtration
generated by an $H$-cylindrical Brownian motion
$(W_H(t))_{t\in [0,T]}$ on a probability space $(\Om,\F,\P)$.
We prove that if $E$ is a UMD Banach space, $1\le p<\infty$,
and $F\in \dD^{1,p}(\Om;E)$ is $\F_T$-measurable, then
$$ F = \E (F) + \int_0^T P_{\mathbb{F}} (DF)\,dW_H,$$
where $D$ is the Malliavin derivative of $F$ and
$P_{\mathbb{F}}$ is the projection onto the
${\mathbb{F}}$-adapted elements in a suitable Banach space of
$L^p$-stochastically integrable $\calL(H,E)$-valued processes.

\end{abstract}

\keywords{Clark-Ocone formula, UMD spaces, Malliavin calculus, Skorokhod integral,
$\g$-radonifying operators, $\g$-boundedness, Stein inequality}
\subjclass[2000]{Primary: 60H07; Secondary: 46B09, 60H05}
\maketitle

\section{Introduction}

A classical result of Clark \cite{Clark} and Ocone
\cite{Ocone} asserts that if $\mathbb{F}=(\F_t)_{t\in [0,T]}$
is the augmented filtration generated by a Brownian motion
$(W(t))_{t\in [0,T]}$ on a probability space $(\Om,\F,\P)$,
then every $\F_T$-measurable random variable $F\in
\dD^{1,p}(\Om)$, $1<p<\infty$, admits a representation
$$ F = \E (F) + \int_0^T \E(D_t F|\F_t)\,dW_t,$$
where $D_t$ is the Malliavin derivative of $F$ at time $t$. An
extension to $\F_T$-measurable random variables $F\in
\dD^{1,1}(\Om)$ was subsequently given by Karatzas, Ocone, and
Li \cite{KOL}. The Clark-Ocone formula is used in mathematical
finance to obtain explicit expressions for hedging strategies.

The aim of this note is to extend the above results to the
infinite-dimensional setting using the theory of stochastic
integration of $\calL(\H,E)$-valued processes with respect to
$\H$-cylindrical Brownian motions, developed recently by
Veraar, Weis, and the second named author \cite{vNVW07}. Here,
$\H$ is a separable Hilbert space and $E$ is a UMD Banach
space (the definition is recalled below).

For this purpose we study the properties of the Malliavin
derivative $D$ of smooth $E$-valued random variables with
respect to an isonormal process $W$ on a separable Hilbert
space $H$. As it turns out, $D$ can be naturally defined as a
closed operator acting from $L^p(\Om;E)$ to
$L^p(\Om;\g(H,E))$, where $\g(H,E)$ is the operator ideal of
$\g$-radonifying operators from a Hilbert space $H$ to $E$.
Via trace duality, the dual object is the divergence operator,
which is a closed operator acting from $L^p(\Om;\g(H,E))$ to
$L^p(\Om;E)$. In the special case where $H = L^2(0,T;\H)$ for
another Hilbert space $\H$, the divergence operator turns out
to be an extension of the UMD-valued stochastic integral of
\cite{vNVW07}.

The first two main results, Theorems \ref{thm:ClarkOcone} and
\ref{thm:ClarkOcone2}, generalize the Clark-Ocone formula for
Hilbert spaces $E$ and exponent $p=2$ as presented in Carmona
and Tehranchi \cite[Theorem 5.3]{CT06} to UMD Banach spaces
and exponents $1<p<\infty$. The extension to $p=1$ is
contained in our Theorem \ref{thm:ClarkOconeL1}.


Extensions of the Clark-Ocone formula to infinite-dimensional
settings different from the one considered here have been
obtained by various authors, among them Mayer-Wolf and Zakai
\cite{MWZa2,MWZa1}, Osswald \cite{Ossw} in the setting of
abstract Wiener spaces and de Faria, Oliveira, Streit
\cite{FOS} and Aase, {\O}ksendal, Privault, Ub{\o}e
\cite{AOPU} in the setting of white noise analysis. Let us
also mention the related papers \cite{M08,MN93}.

\medskip
{\em Acknowledgment} -- Part of this work was done while the
authors visited the University of New South Wales (JM) and the
Australian National University (JvN). They thank Ben Goldys at
UNSW and Alan McIntosh at ANU for their kind hospitality.

\section{Preliminaries}

We begin by recalling some well-known facts concerning
$\g$-radonifying operators and UMD Banach spaces.

Let $(\g_n)_{n\ge 1}$ be sequence of independent standard
Gaussian random variables on a probability space $(\Om,\F,\P)$
and let $H$ be a separable real Hilbert space. A bounded
linear operator $R:H\to E$ is called {\em $\g$-radonifying} if
for some (equivalently, for every) orthonormal basis
$(h_n)_{n\ge 1}$ the Gaussian sum $ \sum_{n\ge 1} \g_n Rh_n$
converges in $L^2(\Om;E)$. Here, $(\g_n)_{n\ge 1}$ is a
sequence of independent standard Gaussian random variables on
$(\Om,\F,\P)$. Endowed with the norm
$$ \n R\n_{\g(H,E)} := \Big(\E\Big\n\sum_{n\ge 1} \g_n Rh_n\big\n^2\Big)^\frac12,
$$
the space $\g(H,E)$ is a Banach space. Clearly $H\otimes
E\subseteq \g(H,E)$, and this inclusion is dense. We have
natural identifications $\g(H,\R) = H$ and $\g(\R,E)=E$.

For all finite rank operators $T:H\to E$ and $S:H\to E\s$ we
have
$$ |{\rm tr} (S\s T)| \le \n T\n_{\g(H,E)}\n
 S\n_{\g(H,E\s)}.
$$
Since the finite rank operators are dense in $\g(H,E)$ and
$\g(H,E\s)$, we obtain a natural contractive injection
\beq\label{eq:traceduality}\g(H,E\s)\hookrightarrow(\g(H,E))\s.\eeq

Let $1<p<\infty$. A Banach space $E$ is called a {\em
UMD$(p)$-space} if there exists a constant $\beta_{p,E}$ such
that for every finite $L^p$-martingale difference sequence
$(d_j)_{j=1}^n$ with values in $E$ and every $\{-1,
1\}$-valued sequence $(\e_j)_{j=1}^n$ we have
$$
\Bigl(\E\Bigl\|\sum_{j=1}^n \e_j d_j \Bigr\|^p\Bigr)^\frac1p
\leq \beta_{p,E} \, \Bigl(\E\Bigl\|\sum_{j=1}^n
d_j\Bigr\|^p\Bigr)^\frac1p.
$$
Using, for instance, Burkholder's good $\l$-inequalities, it
can be shown that if $E$ is a UMD$(p)$ space for some
$1<p<\infty$, then it is a UMD$(p)$-space for all
$1<p<\infty$, and henceforth a space with this property will
simply be called a {\em UMD space}.

Examples of UMD spaces are all Hilbert spaces and the spaces
$L^p(S)$ for $1<p<\infty$ and $\sigma$-finite measure spaces
$(S,\Sigma,\mu)$. If $E$ is a UMD space, then $L^p(S;E)$ is a
UMD space for $1<p<\infty$. For an overview of the theory of
UMD spaces and its applications in vector-valued stochastic
analysis and harmonic analysis we recommend Burkholder's
review article \cite{Bu}.

Below we shall need the fact that if $E$ is a UMD space, then
trace duality establishes an isomorphism of Banach spaces
$$\g(H,E\s)\simeq (\g(H,E))\s.$$
As we shall briefly explain, this is a consequence of the fact
that every UMD is $K$-convex.

Let $(\g_n)_{n\ge 1}$  be sequence of independent standard
Gaussian random variables on a probability space
$(\Om,\F,\P)$. For a random variable $X\in L^2(\Om;E)$ we
define
$$\pi_N^E X:= \sum_{n=1}^N \g_n \E (\g_n X).$$
Each $\pi_N^E$ is a projection on $L^2(\Om;E)$. The Banach
space $E$ is called {\em $K$-convex} if
$$ K(E):= \sup_{N\ge 1} \n \pi_N^E\n <
\infty.$$ In this situation, $\pi^E f :=
\lim_{N\to\infty}\pi_N^E$ defines a projection on $L^2(\Om;E)$
of norm $\n \pi^E\n = K(E)$. It is easy to see that $E$ is
$K$-convex if and only its dual $E\s$ is $K$-convex, in which
case one has $K(E) = K(E\s)$. For more information we refer to
the book of Diestel, Jarchow, Tonge \cite{DJT}.

 The next result from
\cite{Pis89}  (see also \cite{KWsf}) shows that if $E$ is
$K$-convex, the inclusion \eqref{eq:traceduality} is actually
an isomorphism:

\bprop\label{prop:Kconvex-gamma-dual} If $E$ is $K$-convex,
then trace duality establishes an isomorphism of Banach spaces
$$\g(H,E\s)\simeq (\g(H,E))\s.$$
\eprop The main step is to realize that $K$-convexity implies
that the ranges of $\pi^E$ and $\pi^{E\s}$ are canonically
isomorphic as Banach spaces. This isomorphism is then used to
represent elements of $(\g(H,E))\s$ by elements of
$\g(H,E\s)$.

 \brem
Let us comment on the role of the UMD property in this paper.
The UMD property is crucial for two reasons. First, it implies
the $L^p$-boundedness of the vector-valued stochastic
integral. This fact is used at various places (Lemma
\ref{lem:itoduality}, Theorem \ref{thm:skorohodintegral}).
Second, the UMD property is used to obtain the boundedness of
the adapted projection (Lemma \ref{lem:adaptedprojection}).
The results in Sections \ref{sec:Malliavinderivative} and
\ref{sec:divergence} are valid for arbitrary Banach spaces.
 \erem

\section{The Malliavin derivative} \label{sec:Malliavinderivative}

Throughout this note, $(\Om,\F,\P)$ is a complete probability
space, $H$ is a separable real Hilbert space, and $W: H\to
L^2(\Om)$ is an {\em isonormal Gaussian process}, i.e., $W$ is
a bounded linear operator from $H$ to $L^2(\Om)$ such that the
random variables $W(h)$ are centred Gaussian and satisfy
$$ \E (W(h_1)W(h_2)) = [h_1,h_2]_H, \quad h_1, h_2 \in H.$$

A {\em smooth random variable} is a function
$F:\Omega\rightarrow\R$ of the form
 $$F = f(W(h_1), \ldots, W(h_n))$$
with  $f\in C_{\rm b}^{\infty}(\R^n)$ and $h_1,\ldots,h_n\in
H$. Here, $C_{\rm b}^{\infty}(\R^n)$ denotes the vector space
of all bounded real-valued $C^\infty$-functions on $\R^n$
having bounded derivatives of all orders. We say that $F$ is
{\em compactly supported} if $f$ is compactly supported. The
collections of all smooth random variables and compactly
supported smooth random variables are denoted by $\calS$ and
$\calSc$, respectively.

Let $E$ be an arbitrary real Banach space and let $1\le
p<\infty$. Noting that $\calSc$ is dense in $L^p(\Om)$ and
that $L^p(\Om)\otimes E$ is dense in $L^p(\Om;E)$, we see:

\blem $\calSEc$ is dense in $L^p(\Om;E)$. \elem

The {\em Malliavin derivative}  of an $E$-valued smooth random
variable of the form
  $$F = f(W(h_1),\ldots,W(h_n))\otimes x$$ with
 $f\in C_{\rm b}^{\infty}(\R^n)$, $h_1,\ldots,h_n\in H$ and $x\in E$,
 is the random variable $DF:\Om\to \g(H,E)$ defined by
$$
 DF = \sum_{j=1}^n \partial_j f (W(h_1), \ldots,
 W(h_n))\ot (h_j \ot x).
$$
Here, $\partial_j$ denotes the $j$-th partial derivative. The
definition extends to $\calSE$ by linearity.

For $h\in H$ we define $DF(h):\Om\to E$ by $(DF(h))(\om):=
(DF(\om))h$. The following result is the simplest case of the
integration by parts formula. We omit the proof, which is the
same as in the scalar-valued case \cite[Lemma 1.2.1]{Nu06}.

 \blem \label{lem:expDFh}
 For all  $F\in\calSE$ and $h\in H$ we have
 $\E(DF(h)) = \E(W(h)F)$.
 \elem

 A straightforward calculation shows that the following product rule
holds for  $F \in \calSE$ and $G\in \calSEs$:
 \begin{align}\label{eq:prodruleEEstar}
  D\ip{F,G} = \ip{DF,G} + \ip{F,DG}.
 \end{align}
On the left hand side $\ip{\cdot,\cdot}$ denotes the duality
between $E$ and $E^*,$ which is evaluated pointwise on $\Om$.
In the first term on the right hand side, the $H$-valued
pairing $\lb\cdot,\cdot\rb$ between $\g(H,E)$ and $E\s$ is
defined by $\lb R,x\s\rb := R\s x\s$. Similarly, the second
term contains the $H$-valued pairing between $E$ and
$\g(H,E\s),$ which is defined by $\lb x,S\rb := S\s x,$
thereby considering $x$ as an element of $E^{**}.$

For scalar-valued functions $F \in \calS$ we may identify $DF
\in L^2(\Om;\g(H,\R))$ with the classical Malliavin derivative
$DF \in L^2(\Om;H).$ Using this identification we obtain the
following product rule for $F \in \calS$ and $G \in \calSE$:
\begin{align}\label{eq:prodruleER}
 D(FG) = F\, DG + DF \ot G.
\end{align}
An application of Lemma \ref{lem:expDFh} to the product
$\ip{F,G}$ yields the following integration by parts formula
for $F \in \calSE$ and $G \in \calSEs$:
 \begin{align}\label{eq:ibpprel}
 \E \ip{DF(h),G} = \E(W(h)\ip{F,G}) - \E\ip{F,DG(h)}.
 \end{align}
From the identity \eqref{eq:ibpprel} we obtain the following
proposition.

 \bprop \label{prop:closable}
For all $1 \leq p < \infty$, the Malliavin derivative $D$ is
closable as an operator from $L^p(\Om;E)$ into
$L^p(\Om;\g(H,E))$.
 \eprop

 \bpr
Let $(F_n)$ be a sequence in $\calSE$ be such that $F_n \to 0$
in $L^p(\Om;E)$ and $DF_n \to X$ in $L^p(\Om;\g(H,E))$ as
$n\to\infty$. We must prove that $X=0.$

Fix $h \in H$ and define  $$ V_h := \{G \in \calSEs :\ W(h)G
\in \calSEs\}.$$ We claim that $V_h$ is weak$^*$-dense in
$(L^p(\Om;E))^*$. Let $\frac1p+\frac1q=1$. To prove this it
suffices to note that the subspace $\{G\in \calS: W(h) G\in
\calS\}$ is weak$\s$-dense in $L^q(\Om)$ and that $L^q(\Om)\ot
E\s$ is weak$\s$-dense in $(L^p(\Om;E))\s$.

Fix $G\in V_h$. Using \eqref{eq:ibpprel} and the fact that the
mapping $Y \mapsto \E\ip{Y(h), G}$ defines a bounded linear
functional on $L^p(\Om;\g(H,E))$  we obtain
 \begin{align*}
 \E\ip{X(h), G}
 = \lim_{n \to \infty} \E \ip{D F_n(h), G}
 = \lim_{n \to \infty} \E (W(h)\ip{F_n, G}) - \E \ip{F_n, DG(h) }.
 \end{align*}
Since $W(h)G$ and $DG(h)$ are bounded it follows that this
limit equals zero. Since $V_h$ is weak$^*$-dense in
$(L^p(\Om;E))^*$, we obtain that $X(h)$ vanishes almost
surely. Now we choose an orthonormal basis $(h_j)_{j\ge 1}$ of
$H.$ It follows that almost surely we have $X(h_j)=0$ for all
$j\ge 1$. Hence, $X = 0$ almost surely.
 \epr

With a slight abuse of notation we will denote the closure of
$D$ again by $D.$ The domain of this closure in $L^p(\Om;E)$
is denoted by $\dD^{1,p}(\Om;E).$ This is a Banach space
endowed with the norm
 \begin{align*}
 \|F\|_{\dD^{1,p}(\Om;E)} := ( \|F\|_{L^p(\Om;E)}^p
           + \|DF\|_{L^p(\Om;\g(H,E))}^p  )^{\frac1p}.
 \end{align*}
We write $\dD^{1,p}(\Om) := \dD^{1,p}(\Om;\R).$

As an immediate consequence of the closability of the
Malliavin derivative we note that the identities
\eqref{eq:prodruleEEstar}, \eqref{eq:prodruleER},
\eqref{eq:ibpprel} extend to larger classes of functions. This
fact will not be used in the sequel.

 \bprop\label{prop:extended}
Let $1 \le p,q,r < \infty$ such that
$\tfrac{1}{p}+\tfrac{1}{q}=\tfrac{1}{r}.$ \beni
 \item For all $F \in \dD^{1,p}(\Om;E)$ and $G \in \dD^{1,q}(\Om;E^*)$ we
have $\ip{F,G} \in \dD^{1,r}(\Om)$ and
  \begin{align*}
   D\ip{F,G} = \ip{DF,G} + \ip{F,DG}.
  \end{align*}

 \item For all $F \in \dD^{1,p}(\Om)$ and $G \in \dD^{1,q}(\Om;E)$ we
have $FG \in \dD^{1,r}(\Om;E)$ and
 \begin{align*}
 D(FG) = F\, DG + DF \ot G.
 \end{align*}

 \item For all $F \in \dD^{1,p}(\Om;E),$  $G \in \dD^{1,q}(\Om;E^*)$ and $h \in H$
we have $\ip{DF(h),G}\in L^r(\Om)$ and
  \begin{align*}
   \E \ip{DF(h),G} = \E(W(h)\ip{F,G}) - \E\ip{F,D G(h)}.
  \end{align*}
 \een
 \eprop

\section{The divergence operator} \label{sec:divergence}

In this section we construct a vector-valued divergence
operator. The trace inequality \eqref{eq:traceduality} implies
that we have a contractive inclusion  $\g(H,E)
\hookrightarrow(\g(H,E^*))^*.$ Hence for $1<p<\infty$ and
$\tfrac{1}{p}+\tfrac{1}{q}=1,$ we obtain a natural embedding
$$ L^p(\Om;\g(H,E)) \hookrightarrow (L^q(\Om;\g(H,E\s)))^*. $$
For the moment let $D$ denote the Malliavin derivative on
$L^q(\Om;E\s)$, which is a densely defined closed operator
with domain $\dD^{1,q}(\Om;E^*)$ and taking values in
$L^q(\Om;\g(H,E\s))$. The \emph{divergence operator} $\d$ is
the part of the adjoint operator $D\s$ in $L^p(\Om;\g(H,E))$
mapping into $L^p(\Om;E)$. Explicitly, the domain $ \D_p(\d)$
consists of those $X \in L^p(\Om;\g(H,E))$ for which there
exists an $F_X \in L^p(\Om;E)$ such that
$$  \E \lb X, DG\rb = \E \ip{F_X,G}  \ \text{ for all } G \in \dD^{1,q}(\Om;E^*).
$$
The function $F_X$, if it exists, is uniquely determined, and
we define
$$\d(X):= F_X , \quad X \in \D_p(\d).$$
The divergence operator $\d$ is easily seen to be closed, and
the next lemma shows that it is also densely defined.

 \blem \label{lem:divelementary}
We have $\calS\otimes\g(H,E) \subseteq \D_{p}(\d)$ and
$$
\d(f\ot R) = \sum_{j\geq1} W(h_j)f\otimes Rh_j - R(Df), \qquad
f\in \calS, \ R\in \g(H,E).
$$
Here $(h_j)_{j\ge 1}$ denotes an arbitrary orthonormal basis
of $H.$
 \elem

 \bpr
For $f\in \calS$, $R\in \g(H,E)$, and $G \in \calSEs$ we
obtain, using the integration by parts formula
\eqref{eq:ibpprel} (or Proposition \ref{prop:extended}(iii)),
 \begin{align*}
 \E \lb f\ot R, DG\rb
 & = \sum_{j\geq1} \E \ip{f\ot Rh_j, DG(h_j)}
 \\& = \sum_{j\geq1} \E (W(h_j)\ip{f\ot Rh_j, G})
         - \E \ip{[Df,h_j]_H \ot Rh_j, G}
 \\& =  \E \biggip{\sum_{j\geq1}W(h_j) f\ot Rh_j
      - \sum_{j\geq1} [Df, h_j]_H\ot Rh_j,G}
 \\ &=  \E \biggip{\sum_{j\geq1} W(h_j)f\otimes Rh_j - R(Df), G}.
 \end{align*}
The sum $\sum_{j\geq1}W(h_j)f\ot Rh_j$ converges in
$L^{p}(\Om;E).$ This follows from the Kahane-Khintchine
inequalities and the fact that $(W(h_j))_{j\ge 1}$ is a
sequence of independent standard Gaussian variables; note that
the function $f$ is bounded.
 \epr

%

Using an extension of Meyer's inequalities, for UMD spaces $E$
and $1<p<\infty$ it can be shown that $\d$ extends to a
bounded operator from $\dD^{1,p}(\Om;\g(H,E))$ to
$L^p(\Om;E)$. For details we refer to \cite{M08}.

\section{The Skorokhod integral}

We shall now assume that $H = L^2(0,T;\H)$, where $T$ is a
fixed positive real number and $\H$ is a separable real
Hilbert space. We will show that if the Banach space $E$ is a
UMD space, the divergence operator $\d$ is an extension of the
stochastic integral for adapted $\calL(\H,E)$-valued processes
constructed recently in \cite{vNVW07}. Let us start with a
summary of its construction.

Let  $W_\H = (W_\H(t))_{t \in [0,T]}$ be an $\H$-cylindrical
 Brownian motion on $(\Om,\F,\P)$, adapted to a filtration ${\mathbb{F}} =
(\F_t)_{t \in [0,T]}$ satisfying the usual conditions. The
It\^{o} isometry defines an isonormal process $W:
L^2(0,T;\H)\to L^2(\Om)$ by
 \begin{align*}
W(\phi) := \int_0^T \phi\, dW_\H, \qquad \phi \in L^2(0,T;\H).
 \end{align*}

Following \cite{vNVW07} we say that a process $X: (0,T) \times
\Om \to \g(\H,E)$ is an elementary adapted process with
respect to the filtration ${\mathbb{F}}$ if it is of the form
 \begin{align} \label{eq:elementaryadapted}
 X(t,\om) =  \sum_{i = 1}^m  \sum_{j = 1}^n\one_{(t_{i-1},t_i]}(t)
1_{A_{ij}}(\om)\sum_{k=1}^l h_{k}\otimes x_{ijk},
 \end{align}
where $0\le t_0 < \dots < t_n\le T$, the sets $A_{ij}\in
\F_{t_{i-1}}$ are disjoint for each $j$, and
$h_{k},\dots,h_{k}\in \H$ are orthonormal. The stochastic
integral with respect to $W_{\H}$ of such a process is defined
by
 \begin{align*}
 I(X) := \int_0^T X\,dW_\H
     :=  \sum_{i = 1}^m  \sum_{j =1}^n \sum_{k=1}^l
     1_{A_{ij}}(W_\H(t_i)h_{k} - W_\H(t_{i-1})h_{k}) \otimes x_{ijk},
 \end{align*}

Elementary adapted processes define elements of
$L^p(\Om;\g(L^2(0,T;\H),E))$ in a natural way. The closure of
these elements in $L^p(\Om;\g(L^2(0,T;\H),E))$ is denoted by
$L^p_{{{\mathbb{F}}}}(\Om;\g(L^2(0,T;\H),E)).$

 \bprop[{\cite[Theorem 3.5]{vNVW07}}] \label{prop:vnvw} Let $E$ be a UMD space and let
$1<p<\infty.$ The stochastic integral uniquely extends to a
bounded operator $$I :
L_{{{\mathbb{F}}}}^p(\Om;\g(L^2(0,T;\H),E)) \to L^p(\Om;E).$$
Moreover, for all $X\in
L_{{{\mathbb{F}}}}^p(\Om;\g(L^2(0,T;\H),E))$ we have the
two-sided estimate
 \begin{align*}
 \|I(X)\|_{L^p(\Om;E)} \eqsim
  \|X\|_{L^p(\Om;\g(L^2(0,T;\H),E))},
 \end{align*}
with constants only depending on $p$ and $E$.
 \eprop

A consequence of this result is the following lemma, which
will be useful in the proof of Theorem \ref{thm:ClarkOcone}.

 \blem \label{lem:itoduality}
Let $E$ be a UMD space and let $1 < p,q < \infty$ satisfy
$\frac1p+\frac1q=1$. For all $X \in
L^p_{{{\mathbb{F}}}}(\Om;\g(L^2(0,T;\H),E))$ and $Y \in
L^q_{{{\mathbb{F}}}}(\Om;\g(L^2(0,T;\H),E^*))$ we have
 \begin{align*}
 \E\ip{I(X),I(Y)}  = \E\lb X,Y\rb.
 \end{align*}
 \elem

  \bpr
When $X$ and $Y$ are elementary adapted the result follows by
direct computation. The general case follows from Proposition
\ref{prop:vnvw} applied to $E$ and $E^*,$ noting that $E^*$ is
a UMD space as well.
  \epr

In the next approximation result we identify $L^2(0,t;\H)$
with a closed subspace of $L^2(0,T;\H)$. The simple proof is
left to the reader.

 \blem \label{lem:approxadapted}
Let $1\leq p<\infty,$ let $0 < t \leq T,$ and let
$(\psi_n)_{n\ge 1}$ be an orthonormal basis of $L^2(0,t;\H)$.
The linear span of the functions $f(W(\psi_1), \ldots,
W(\psi_n)) \ot (h\ot x)$, with $f \in \calS$, $h\in H$, $x \in
E$, is dense in $L^p(\Om,\F_t;\g(\H,E))$.
 \elem

The next result shows that the divergence operator $\d$ is an
extension of the stochastic integral $I.$ This means that $\d$
is a vector-valued Skorokhod integral.

 \bthm \label{thm:skorohodintegral}
Let $E$ be a UMD space and let $1<p<\infty$ be fixed. The
space $L^p_{{{\mathbb{F}}}}(\Om;\g(L^2(0,T;\H),E))$ is
contained  in $\D_p(\d)$ and
$$\d(X) = I(X) \ \hbox{ for all } \ X \in L^p_{{{\mathbb{F}}}}(\Om;\g(L^2(0,T;\H),E)).$$
 \ethm

 \bpr
Fix  $0 < t \leq T$, let $(h_k)_{k\ge 1}$  be an orthonormal
basis of $\H$, and put $X := 1_A \sum_{k=1}^n h_k\ot x_k$ with
$A\in \F_{t}$ and $x_k \in E$ for $k=1,\dots, n$. Let
$(\psi_j)_{j\ge 1}$ be an orthonormal basis of $L^2(0,t;\H)$.
By Lemma \ref{lem:approxadapted} we can approximate $X$ in
$L^p(\Om,\F_t;\g(\H,E))$ with a sequence $(X_l)_{l\ge 1}$ in
$\mathscr{S}(\Om,\g(\H,E))$ of the form
 \begin{align*}
 X_l := \sum_{m=1}^{M_l} f_{lm}(W(\psi_{1}), \ldots , W(\psi_{n})) \ot
 (h_{m}\ot x_{lm})
 \end{align*}
with $x_{lm}\in E$.

Now let $0<t<u\le T$. From $\psi_{m}\perp \one_{(t,u]}\ot h$
in $L^2(0,T;\H)$ it follows that  $ DX_l(\one_{(t,u]}\ot h) =
0$ for all $h\in\H$. By Lemma \ref{lem:divelementary},
$$\one_{(t,u]} \ot X_l =
\sum_{m=1}^{M_l} f_{lm}(W(\psi_{1}), \ldots , W(\psi_{n})) \ot
((\one_{(t,u]}\ot h_{m})\ot x_{lm})
$$
belongs to $\D_p(\d)$
and 
$$
\d(\one_{(t,u]} \ot X_l )
  = \sum_{m=1}^{M_l} W(\one_{(t,u]}\ot h_{m})
f_{lm}(W(\psi_{1}), \ldots , W(\psi_{n})) \ot x_{lm}
 = I(\one_{(t,u]}\ot X_l).
$$
Noting that $ \one_{(t,u]} \ot X_l \to \one_{(t,u]} \ot X$ in
$L^p(\Om;\g(L^2(0,T;\H),E))$ as $l\to\infty$, the closedness
of $\d$ implies that $\one_{(t,u]} \ot X \in \D_p(\d)$ and
$$ \d( \one_{(t,u]} \ot X_l ) = I(\one_{(t,u]}\ot X_l).$$
By linearity, it follows that the elementary adapted processes
of the form \eqref{eq:elementaryadapted} with $t_0>0$ are
contained in $\D_p(\d)$ and that $I$ and $\d$ coincide for
such processes.

To show that this equality extends to all $ X \in
L^p_{{{\mathbb{F}}}}(\Om;\g(L^2(0,T;\H),E))$ we take a
sequence $X_n$ of elementary adapted processes of the above
form converging to $X.$ Since $I$ is a bounded operator from
$L^p_{{{\mathbb{F}}}}(\Om;\g(L^2(0,T;\H),E))$ into
$L^p(\Om;E),$ it follows that $\d(X_n) = I(X_n) \to I(X)$ as
$n\to\infty$. The fact that $\d$ is closed implies that $X \in
\D_p(\d)$ and $\d(X) = I(X).$
 \epr

\section{A Clark-Ocone formula}\label{sec:Clark-Ocone}

Our next aim is to prove that the space
$L_{{{\mathbb{F}}}}^p(\Om;\g(L^2(0,T;\H),E))$, which has been
introduced in the previous section, is complemented in
$L^p(\Om;\g(L^2(0,T;\H),E))$. For this we need a number of
auxiliary results. Before we can state these we need to
introduce some terminology.

Let $(\gamma_j)_{j\ge 1}$ be a sequence of independent
standard Gaussian random variables. Recall that a collection
$\mathscr{T} \subseteq \L(E,F)$ of bounded linear operators
between Banach spaces $E$ and $F$ is said to be {\em
$\g$-bounded} if there exists a constant $C>0$ such that
 \begin{align*}
 \E \Big\| \sum_{j=1}^n \g_j T_j x_j \Big\|_F^2
  \le C^2\E \Big\| \sum_{j=1}^n \g_j  x_j \Big\|_E^2
 \end{align*}
for all $n \ge 1$ and all choices of $T_1,\ldots,T_n \in
\mathscr{T}$ and $x_1,\ldots, x_n \in E.$ The least admissible
constant $C$ is called the {\em $\g$-bound} of $\mathscr{T}$,
notation $\g(\mathscr{T})$.

\bprop\label{prop:lifting} Let $\mathscr{T}$ be a $\g$-bounded
subset of $\calL(E,F)$ and let $H$ be a separable real Hilbert
space. For each $T\in \mathscr{T}$ let $\widetilde T\in
\calL(\g(H,E),\g(H,F))$ be defined by $\widetilde TR := T\circ
R$. The collection $\widetilde{\mathscr{T}} = \{\widetilde T:\
T\in \mathscr{T}\}$ is $\g$-bounded, with
$\g(\widetilde{\mathscr{T}}) = \g(\mathscr{T})$. \eprop

\begin{proof}
Let $(\gamma_j)_{j\ge 1}$ and $(\widetilde\gamma_j)_{j\ge 1}$
be two sequences of independent standard Gaussian random
variables, on probability spaces $(\Om,\F,\P)$ and
$(\widetilde\Om, \widetilde\F,\widetilde \P)$ respectively. By
the Fubini theorem,
$$
\begin{aligned}
\E \Big\n \sum_{j=1}^n \g_j \widetilde T_j R_j
\Big\n_{\g(H,F)}^2 & = \E \widetilde{\E} \Big\n
\sum_{i=1}^\infty \widetilde\g_i \sum_{j=1}^n \g_j T_j R_j h_i
\Big\n_F^2
\\ & =  \widetilde\E \E \Big\n \sum_{j=1}^n \g_j T_j\sum_{i=1}^\infty \widetilde\g_i
R_j h_i \Big\n_F^2
 \\ & \le \g^2(\mathscr{T})\widetilde \E\E\Big\n \sum_{j=1}^n\g_j \sum_{i=1}^\infty
 \widetilde \g_i R_j h_i \Big\n_E^2
\\ & = \g^2(\mathscr{T}) \E \widetilde \E \Big\n \sum_{i=1}^\infty \widetilde
\g_i \sum_{j=1}^n \g_j R_j h_i \Big\n_E^2
\\ & = \g^2(\mathscr{T})\E \Big\n \sum_{j=1}^n \g_j R_j \Big\n_{\g(H,E)}^2.
\end{aligned}
$$
This proves the inequality $\g(\widetilde{\mathscr{T}}) \le
\g(\mathscr{T})$. The reverse inequality holds trivially.
\end{proof}

The next proposition is a result by Bourgain \cite{Bo86},
known as the vector-valued Stein inequality. We refer to
\cite[Proposition 3.8]{CPSW00} for a detailed proof.

 \bprop \label{prop:Steinineq}
Let $E$ be a UMD space and let $(\F_t)_{t \in [0,T]}$ be a
filtration on $(\Om,\F,\P)$. For all $1 < p < \infty$  the
conditional expectations $\{ \E(\cdot|\F_t) :\ t \in [0,T]\}$
define a $\g$-bounded set in $\calL(L^p(\Om;E))$.
 \eprop

We continue with a multiplier result due to Kalton and Weis
\cite{KWsf}. In its formulation we make the observation that
every step function $f:(0,T)\to \g(\H,E)$ defines an element
$R_g\in \g(L^2(0,T;\H),E)$ by the formula
$$ R_f \phi := \int_0^T f(t)\phi(t)\,dt.$$
Since $R_f$ determines $f$ uniquely almost everywhere, in what
follows we shall always identify $R_f$ and $f$.

 \bprop \label{prop:KaltonWeismultiplier}
Let $E$ and $F$ be real Banach spaces and let $M: (0,T) \to
\calL(E,F)$ have $\g$-bounded range $\{M(t) : t \in (0,T)\}
=:\mathscr{M}$. Assume that for all $x\in E$, $t\mapsto M(t)x$
is strongly measurable. Then the mapping $M: f\mapsto
[t\mapsto M(t)f(t)]$ extends to a bounded operator from
$\g(L^2(0,T;\H),E)$ to $\g(L^2(0,T;\H),F)$ of norm $\n M\n\le
\gamma(\mathscr{M})$.
 \eprop
Here we identified $M(t)\in\calL(E,F)$ with
$\widetilde{M(t)}\in\calL(\g(\H,E),\g(\H,F))$ as in
 Proposition \ref{prop:lifting}.

The next result is taken from \cite{vNVW07}.

 \bprop \label{prop:gLpisomorphicLpg}
Let $H$ be a separable real Hilbert space and let $1\leq p <
\infty$. Then $ f\mapsto [h\mapsto f(\cdot)h]$ defines an
iso\-morphism of Banach spaces $$L^p(\Om;\g(H,E))\simeq
\g(H,L^p(\Om;E)).$$ \eprop

After these preparations we are ready to state the result
announced above. We fix a filtration  $\mathbb{F} = (\F_t)_{t
\in [0,T]}$ and define, for step functions $f: (0,T)\to
\g(\H,L^p(\Om;E))$, \beq\label{eq:defPF}
 (P_{\mathbb{F}} f)(t) :=\E(f(t) | \F_t),
\eeq where $\E(\cdot|\F_t)$ is considered as a bounded
operator acting on $\g(\H,L^p(\Om;E))$ as in Proposition
\ref{prop:lifting}.

\blem \label{lem:adaptedprojection} Let $E$ be a UMD space,
and let $1 < p,q < \infty$ satisfy $\frac1p+\frac1q = 1.$
 \beni
 \item
The mapping $ P_{{{\mathbb{F}}}}$ extends to a bounded
operator on $\g(L^2(0,T;\H),L^p(\Om;E))$.
 \item
As a bounded operator on $L^p(\Om;\g(L^2(0,T;\H),E))$,
$P_{{{\mathbb{F}}}}$ is a projection onto the subspace
$L_{{{\mathbb{F}}}}^p(\Om;\g(L^2(0,T;\H),E)).$
 \item
For all $X \in L^p(\Om;\g(L^2(0,T;\H),E))$ and $Y \in
L^q(\Om;\g(L^2(0,T;\H),E^*))$ we have $$\E\lb X,
P_{{{\mathbb{F}}}}Y\rb = \E\lb P_{{{\mathbb{F}}}}X, Y\rb.$$
 \item
For all $X \in L^p(\Om;\g(L^2(0,T;\H),E))$ we have $\E
P_{\mathbb{F}}X=\E X$.
 \een
 \elem

 \bpr
(i): From Propositions \ref{prop:lifting} and
\ref{prop:Steinineq} we infer that the collection of
conditional expectations $\{\E(\cdot|\F_t) : t \in [0,T]\}$ is
$\g$-bounded in $\L(\g(\H,L^p(\Om;E)))$. The boundedness of
$P_{{{\mathbb{F}}}}$ then follows from Proposition
\ref{prop:KaltonWeismultiplier}. For step functions
$f:(0,T)\to \g(\H,L^p(\Om;E))$ it is clear from
\eqref{eq:defPF} that $P_{{{\mathbb{F}}}}^2 f =
P_{{{\mathbb{F}}}}f$, which means that $P_{{{\mathbb{F}}}}$ is
a projection.\\ \smallskip

(ii): By the identification of Proposition
\ref{prop:gLpisomorphicLpg}, $P_{\mathbb{F}}$ acts as a
bounded projection in the space $L^p(\Om;\g(L^2(0,T;\H),E))$.
For elementary adapted processes $X\in
L^p(\Om;\g(L^2(0,T;\H),E))$ we have $P_{\mathbb{F}}X=X$, which
implies that the range of $P_{{{\mathbb{F}}}}$ contains $
L_{{{\mathbb{F}}}}^p(\Om;\g(L^2(0,T;\H),E))$. To prove the
converse inclusion we fix a step function $X:(0,T)\to \g(\H,
L^p(\Om;E))$ and observe that $P_{{{\mathbb{F}}}}X$ is adapted
in the sense that $(P_{{{\mathbb{F}}}}X)(t)$ is strongly
$\F_t$-measurable for every $t\in [0,T].$ As is shown in
\cite[Proposition 2.12]{vNVW07}, this implies that
$P_{{{\mathbb{F}}}} X \in
L_{{{\mathbb{F}}}}^p(\Om;\g(L^2(0,T;\H),E)).$ By density it
follows that the range of $P_{\mathbb{F}}$ is contained in
$L_{{{\mathbb{F}}}}^p(\Om;\g(L^2(0,T;\H),E)).$ \\ \smallskip

(iii): Keeping in mind the identification of Proposition
\ref{prop:gLpisomorphicLpg}, for step functions with values in
the finite rank operators from $\H$ to $E$ this follows from
\eqref{eq:defPF} by elementary computation. The result then
follows from a density argument. \\ \smallskip

(iv): Identifying a step function $f:(0,T)\to
\g(\H,L^p(\Om;E))$ with the associated operator in
$\g(L^2(0,T;\H),L^p(\Om;E))$ and viewing $\E$ as a bounded
operator from $\g(L^2(0,T;\H),L^p(\Om;E))$ to
$\g(L^2(0,T;\H),E)$, by \eqref{eq:defPF} we have
$$ \E P_{\mathbb{F}}f(t) = \E \E(f(t)|\F_t) = \E f(t).$$
Thus $\E P_{\mathbb{F}}f = \E f$ for all step functions
$f:(0,T)\to\g(\H,L^p(\Om;E))$, and hence for all $f\in
\g(L^2(0,T;\H),L^p(\Om;E))$ by density. The result now follows
by an application of Proposition \ref{prop:gLpisomorphicLpg}.
\epr

Now let $\mathbb{F} = (\F_t)_{t\in [0,T]}$ be the augmented
filtration generated by $W_{\H}$. It has been proved in
\cite[Theorem 4.7]{vNVW07} that if  $E$ is a UMD space and $1
< p < \infty,$ and if $F\in L^p(\Om;E)$ is $\F_T$-measurable,
then there exists a unique $X \in
L_{{{\mathbb{F}}}}^p(\Om;\g(L^2(0,T;\H),E))$ such that
$$F = \E(F) + I(X).$$
The following two results give an explicit expression for $X$.
They extend the classical Clark-Ocone formula and its Hilbert
space extension to UMD spaces.

 \bthm[Clark-Ocone representation, first $L^p$-version]\label{thm:ClarkOcone}
Let $E$ be a UMD space and let $1 < p < \infty.$ If $F \in
\dD^{1,p}(\Om;E)$ is $\F_T$-measurable, then
 \begin{align*}
 F = \E(F) + I (P_{\mathbb{F}}(D F)).
 \end{align*}
Moreover,  $P_{\mathbb{F}}(D F)$ is the unique $Y\in
L_{\mathbb{F}}^p(\Om;\g(L^2(0,T;\H),E))$ satisfying $F = \E(F)
+ I(Y)$.
 \ethm

\begin{proof}
We may assume that $\E(F) = 0.$ Let  $X\in
L_{{{\mathbb{F}}}}^p(\Om;\g(L^2(0,T;\H),E))$ be such that $F =
I(X) = \d(X).$ Let $\frac1p+\frac1q= 1,$ and let $Y \in
L^q(\Om;\g(L^2(0,T;\H),E^*))$ be arbitrary. By Lemma
\ref{lem:adaptedprojection}, Theorem
\ref{thm:skorohodintegral}, and Lemma \ref{lem:itoduality} we
obtain
 \begin{align*}
    \E\lb P_{{{\mathbb{F}}}}(DF),Y\rb
& =  \E \lb DF,P_{{{\mathbb{F}}}} Y\rb
 =  \E\lb F,\d(P_{{{\mathbb{F}}}} Y)\rb
\\ & =  \E\lb\d(X),\d(P_{{{\mathbb{F}}}} Y)\rb
 =  \E\lb I(X),I(P_{{{\mathbb{F}}}} Y)\rb
\\& =  \E \lb X,P_{{{\mathbb{F}}}} Y\rb
 =  \E\lb P_{{{\mathbb{F}}}}X, Y\rb
 =  \E\lb X, Y\rb.
 \end{align*}
Since this holds for all $Y \in L^q(\Om;\g(L^2(0,T;\H),E^*)),$
it follows that $X = P_{{{\mathbb{F}}}}(DF).$ The uniqueness
of $P_{{{\mathbb{F}}}}(DF)$ follows from the injectivity of
$I$ as a bounded linear operator from
$L_{{{\mathbb{F}}}}^p(\Om;\g(L^2(0,T;\H),E))$ to
$L^p(\Om,\F_T)$.
\end{proof}

With a little extra effort we can prove a bit more:

 \bthm[Clark-Ocone representation, second $L^p$-version]\label{thm:ClarkOcone2}
Let $E$ be a UMD space and let $1 < p < \infty.$ The operator
$P_{\mathbb F}\circ D$ has a unique extension to a bounded
operator from $L^p(\Om,\F_T;E)$ to
$L_{{{\mathbb{F}}}}^p(\Om;\g(L^2(0,T;\H),E))$, and for all
$F\in L^p(\Om,\F_T;E)$ we have the representation
 \begin{align*}
 F = \E(F) + I ((P_{{{\mathbb{F}}}}\circ D)F).
 \end{align*}
Moreover,  $(P_{\mathbb F}\circ D)F$ is the unique $Y\in
L_{\mathbb{F}}^p(\Om;\g(L^2(0,T;\H),E))$ satisfying $F = \E(F)
+ I(Y)$.
 \ethm
\begin{proof}
It follows from Theorem \ref{thm:ClarkOcone} that  $F\mapsto
I((P_{{{\mathbb{F}}}}\circ D)F) $ extends uniquely to a
bounded operator on $L^p(\Om,\F_T;E)$, since it equals
$F\mapsto F - \E(F)$ on the dense subspace $
\dD^{1,p}(\Om,\F_T;E)$. The proof is finished by recalling
that $I$ is an isomorphism from
$L_{{{\mathbb{F}}}}^p(\Om;\g(L^2(0,T;\H),E))$ onto its range
in $L^p(\Om,\F_T)$.
\end{proof}

\brem An extension of the Clark-Ocone formula to a class of
adapted processes taking values in an arbitrary Banach space
$B$ has been obtained by Mayer-Wolf and Zakai \cite[Theorem
3.4]{MWZa2}. The setting of \cite{MWZa2} is slightly different
from ours in that the starting point is an arbitrary abstract
Wiener space $(W,H,\mu)$, where $\mu$ is a centred Gaussian
Radon measure on the Banach space $W$ and $H$ is its
reproducing kernel Hilbert space. The filtration is defined in
terms of an increasing resolution of the identity on $H$, and
a somewhat weaker notion of adaptedness is used. However, the
construction of the predictable projection in \cite[Section
3]{MWZa2} as well as the proofs of \cite[Corollary 3.5 and
Proposition 3.14]{MWZa1} contain gaps. As a consequence, the
Clark-Ocone formula of \cite{MWZa2} only holds in a suitable
`scalar' sense. We refer to the errata \cite{MWZa2, MWZa1} for
more details.
 \erem

\section{Extension to $L^1$}

We continue with an extension of Theorem \ref{thm:ClarkOcone2}
to random variables in the space $L^1(\Om,\F_T;E)$. As before,
$\mathbb{F} = (\F_t)_{t\in [0,T]}$ is the augmented filtration
generated by the $\H$-cylindrical Brownian motion $W_\H$.

We denote by $L^0(\Om;F)$ the vector space of all strongly
measurable random variables with values in the Banach space
$F$, identifying random variables that are equal almost
surely. Endowed with the metric
$$d(X,Y) = \E(\n X-Y\n\wedge 1),$$ $L^0(\Om;F)$ is a complete metric space,
and we have $\lim_{n\to\infty} X_n=X$ in $L^0(\Om;F)$ if and
only if $\lim_{n\to\infty} X_n = X$ in measure in $F$.

The closure of the elementary adapted processes in
$L^0(\Om;\g(L^2(0,T;\H),E))$ is denoted by
$L^0_{\mathbb{F}}(\Om;\g(L^2(0,T;\H),E)).$ By the results of
\cite{vNVW07}, the stochastic integral $I$ has a unique
extension to a linear homeomorphism from
$L_{\mathbb{F}}^0(\Om;\g(L^2(0,T;\H),E))$ onto its image in
$L^0(\Om,\F_T;E)$.

 \bthm[Clark-Ocone representation, $L^1$-version]\label{thm:ClarkOconeL1}
Let $E$ be a UMD space. The operator $P_{\mathbb F}\circ D$
has a unique extension to a continuous linear operator from
$L^1(\Om,\F_T;E)$ to
$L_{{{\mathbb{F}}}}^0(\Om;\g(L^2(0,T;\H),E))$, and for all
$F\in L^1(\Om,\F_T;E)$ we have the representation
 \begin{align*}
 F = \E(F) + I ((P_{{{\mathbb{F}}}}\circ D)F).
 \end{align*}
Moreover,  $(P_{\mathbb F}\circ D)F$ is the unique element
$Y\in L_{\mathbb{F}}^0(\Om;\g(L^2(0,T;\H),E))$ satisfying $F =
\E(F) + I(Y)$.
 \ethm

\begin{proof}
We shall employ the process $\xi_X:(0,T)\times\Om\to
\g(L^2(0,T;\H),E)$ associated with a strongly measurable
random variable $X: \Om \to \g(L^2(0,T;\H),E)$, defined by
$$(\xi_X(t,\om))f:= (X(\om))(\one_{[0,t]}f), \qquad f \in
L^2(0,T;\H).$$ Some properties of this process have been
studied in  \cite[Section 4]{vNVW07}.

Let $(F_n)_{n\ge 1}$ be a sequence of $\F_T$-measurable random
variables in $\mathscr{S}(\Om)\otimes E$ which is Cauchy in
$L^{1}(\Om,\F_T;E)$. By \cite[Lemma 5.4]{vNVW07} there exists
a constant $C\ge 0$, depending only on $E$, such that for all
$\d>0$ and $\e>0$ and all $m,n\ge 1$,
$$
\begin{aligned}
\ &  \P\big(\n P_{\mathbb{F}}
(DF_n-DF_m)\n_{\g(L^2(0,T;\H),E)} > \e \big)
\\ & \qquad\qquad \le \frac{C\d^2}{\e^2}
+ \P\big(\sup_{t\in [0,T]}\n I(\xi_{P_{\mathbb{F}}
(DF_n-DF_m)}(t))\n \ge \d \big)
\\ &\qquad\qquad \stackrel{(*)}{=} \frac{C\d^2}{\e^2}
+ \P\big(\sup_{t\in [0,T]}\n \E(F_n- F_m|\F_t)-\E(F_n- F_m)\n
\ge  \d \big)
\\ &\qquad\qquad \stackrel{(**)}{\le}  \frac{C\d^2}{\e^2}
+ \frac1\d \E \n F_n- F_m-\E(F_n- F_m)\n.
\end{aligned}
$$
In this computation, $(*)$
 follows from Theorem \ref{thm:ClarkOcone} which gives
$$ \E(F|\F_t) - \E (F)
   =   \E \big(I (P_{\mathbb{F}}DF)\big|\F_t\big)
   =  \E \big(I (\xi_{P_{\mathbb{F}}DF}(T))\big|\F_t\big)
   = I (\xi_{P_{\mathbb{F}}DF}(t)).$$
The estimate $(**)$ follows from Doob's maximal inequality.
Since the right-hand side in the above computation can be made
arbitrarily small, this proves that
$(P_{\mathbb{F}}(DF_n))_{n\ge 1}$ is Cauchy in measure in
$\g(L^2(0,T;\H),E)$.

For $F\in L^1(\Om,\F_T;E)$ this permits us to define
$$ (P_{\mathbb{F}} \circ D)F:= \lim_{n\to\infty} P_{\mathbb{F}}(DF_n),$$
where $(F_n)_{n\ge 1}$ is any sequence of $\F_T$-measurable
random variables in $\calSE$ satisfying $\lim_{n\to\infty} F_n
= F$ in $L^{1}(\Om,\F_T;E)$. It is easily checked that this
definition is independent of the approximation sequence. The
resulting linear operator $P_{\mathbb{F}} \circ D$ has the
stated properties. This time we use the fact that $I$ is a
homeomorphism from $L_{\mathbb{F}}^0(\Om;\g(L^2(0,T;\H),E))$
onto its image in $L^0(\Om,\F_T;E)$; this also gives the
uniqueness of $(P_{\mathbb{F}} \circ D)F$.
\end{proof}

\bibliographystyle{ams-pln}
\bibliography{ecp_clark_ocone}

\end{document}